\documentclass[12pt]{amsart}
\usepackage{amsmath}
\usepackage{amsthm}
\usepackage{amssymb}

\newtheorem{theorem}{Theorem}

\theoremstyle{definition}

\theoremstyle{remark}

\numberwithin{equation}{section}

\begin{document}

\title{Additive properties of even perfect numbers}
\author{Yu Tsumura}
\date{}

\begin{abstract}
A positive integer $n$ is said to be perfect if $\sigma(n)=2n$, where $\sigma$ denotes the sum of the divisors of $n$.
In this article, we  show that if $n$ is an even perfect number, then any integer $m\leq n$ is expressed as a sum of some of divisors of $n$.
\end{abstract}
\maketitle
\section{Introduction.}

A perfect number is a positive integer whose proper divisors sum up to $n$ itself.
Euclid proved that if $2^p-1$ is prime, then $n=2^{p-1}(2^p-1)$ is an even perfect number.
Two millennia later, Euler proved the converse.
\begin{theorem}
An even positive integer $n$ is perfect if and only if it is of the form $n=2^{p-1}(2^p-1)$, where $2^p-1$ is prime.
\end{theorem}\label{thm1}
For a proof, see \cite[Theorem 1.3.3, p.\ 22]{Crandall-Pomerance}.
So we have a  characterization of even perfect numbers.
However, it is not known whether there are infinitely many perfect numbers.
In addition, we do not know whether there exists an odd perfect number.
See \cite{Guy} for these problems and related unsolved problems.
\section{Main results.}

If $n$ is a perfect number, then the sum of all its proper divisor is $n$.
Then, what can we say about the sum of some of its proper divisors?
Using Theorem \ref{thm1}, we can prove the following theorem.
\begin{theorem}\label{thm:main}
If $n$ is an even perfect number, then any positive integer $m$ less than or equal to $n$ is expressed as a sum of some of divisors of $n$.
\end{theorem}

Before proving it, we give an example.
$6$ is an even perfect number since its proper divisors $1$, $2$, and $3$ sum up to $6$.
Theorem \ref{thm:main} says every positive integer $m\leq 6$ can be written as a sum of some of $1$, $2$, and $3$.
For $m=1$, $2$, and $3$, they are themselves proper divisors, so we do not need to think about them.
For $m=4$, we have $4=1+3$.
For $m=5$, we have $5=2+3$.
Finally, since $6$ itself is a perfect number, we have $6=1+2+3$.

Let us move on to the proof of  Theorem \ref{thm:main}

\begin{proof}
From now on, any divisor means a proper divisor of $n$.
Suppose $n$ is an even perfect number.
By Theorem \ref{thm1}, we can write $n=2^{p-1}M_p$, where $M_p=2^p-1$ is prime.
First of all, note that every number $m$ such that $1\leq m \leq 2^p-1$ is a sum of some of $1$, $2$, $2^2$, $\ldots$, $2^{p-1}$. 
(This is just a binary representation of $m$.)
Since $1$, $2$, $2^2$, $\ldots$, $2^{p-1}$ are  divisors of $n$, we could express every number in $S_0=\{1$, $2$, $3$, $\ldots$, $M_p \}$ as a sum of some of  divisors.
Since we did not use a  divisor $M_p$, we can add it to numbers in $S_1$ and we see that every number in $S_1=\{1+M_p$, $2+M_p$, $3+M_p$, $\ldots$, $2M_p\}$ is a sum of some of its divisors.
Next, adding $2M_p$ to numbers in $S_0$, we see that every number in  $S_2= \{1+2M_p$, $2+2M_p$, $3+2M_p$, $\ldots$, $3M_p\}$ is a some of some of its divisors.
Similarly, adding $kM_p$ with $1\leq k \leq 2^{p-1}-1$ to numbers in $S_0$, we see that every number in $S_k=\{1+kM_p$, $2+kM_p$, $3+kM_p$, $\ldots$, $(k+1)M_p\}$ is a sum of some of its divisors.
Since the set of positive integer less than or equal to $n$ is $\bigcup_{k=0}^{2^{p-1}-1}S_k$, we have expressed every $m$ less than or equal to $n$ as a sum of some of divisors of $n$.
\end{proof}

The next question is whether this expression is unique or not.
The answer is negative.
For example, we take a perfect number $n=6$ and $m=3$.
Then we can express $m=1+2=3$ with proper divisors of $6$.
Hence $m=3$ can be expressed in two ways.
In general, since $M_p=\sum_{i=0}^{p-1}2^i$, $M_p$ is expressed in two ways.
Hence a multiple of $M_p$ is expressed in two ways except for $n$ itself.
However, we show that this is only the case in the following theorem.

\begin{theorem}
The expression in Theorem \ref{thm:main} is unique except for $m=kM_p$, where $k=1$, $2$, $3$, $\ldots$, $2^{p-1}-1$.
Also $m=kM_p$ is expressed in exactly two ways for $k=1$, $2$, $3$, $\ldots$, $2^{p-1}-1$.
\end{theorem}
\begin{proof}
First, we have seen that $kM_p$ is expressed in (at least) two ways for $k=1$, $2$, $3$, $\ldots$, $2^{p-1}-1$.
Now counting multiplicity, there are $2^{2p-1}-1$ combinations of divisors of $n$ since the number of proper divisors of $n=2^{p-1}M_p$ is $2p-1$.
(To exclude an empty combination, we subtracted $1$ from $2^{2p-1}$.)
Since $2^{2p-1}-1-n=2^{p-1}-1$, there are no $m$ other than $m=kM_p$ that is expressed in two ways and $kM_p$ is not expressed in more than two ways.
\end{proof}

The next natural question is whether the converse is true.
That is, if $n$ is an even positive integer and if every positive integer $m$ less than or equal to $n$ is expressed as a sum of some of proper divisors of $n$, can we say that $n$ is an even perfect number?
Again, the answer is negative.
we give a counterexample.
Let us consider $n=20$.
Then it is easy to check that $20$ is not a perfect number and some of its proper divisors $1$, $2$, $4$, $5$, $10$ sum up to all positive integers less than or equal to $20$.

In conclusion, divisors of an even perfect number $n$ construct not only $n$ itself but also all positive integers less than or equal to $n$.
What a perfect number it is!

\bibliographystyle{plain}
\bibliography{perfectnumber}

\bigskip

\noindent\textit{Department of Mathematics, Purdue University
150 N. University Street, West Lafayette, IN 47907-2067\\
ytsumura@math.purdue.edu}

\end{document}